\documentclass[12pt]{article}

\title{Minimum settling time control design through direct search optimization}
\author{Emile Simon}		

\usepackage{animate}
\usepackage{graphicx}           

\begin{document}
\maketitle						

The aim of this very brief third version is to present succinctly, but clearly with two animations, two optimization results from the second version of the paper available on http://arxiv.org/abs/1109.5966v2 which can still be consulted for more details. A future and much improved version of the paper will be deposited on arXiv at the beginning of 2012.\\

The objective considered here is to optimize PID coefficients $x=[K_p,K_i,K_d]$ to minimize the following objective function:\\

$f(x)$ = the rise time (= time needed for the response $z(x,t)$ to reach $0.98$) divided by $100$ (the fixed maximum time value for which the response is evaluated)

+ the maximum deviation of $z(x,t)$ (out of the following settling range: $z_{max}(t)=1.02$ for $t>0$ and $z_{min}(t)=0.98$ for $t>$ rise time).

\newpage

First figure of unit step response optimization, starting from PID coefficients obtained by Ziegler-Nichols:

\begin{center}
\animategraphics[controls,loop,width=4in]{12}{Film7Bis/film_}{1}{217}
\end{center}

Press play to start the animation.

NB/ The response $z(x,t)$ is drawn in green when the corresponding objective value $f(x)$ is better than best value previously found by the direct search optimization method, and in red otherwise. The horizontal black dashed lines represent the settling range defined above.

\newpage

Second figure of unit step response optimization, starting from unstabilizing random PID coefficients: 

\begin{center}
\animategraphics[controls,width=4in]{12}{Film8/film_}{1}{437}

time [s]
\end{center}

Both optimziations lead to the same solution, with rise time around 11.94s and maximum deviation almost zero (can be fine-tuned with another objective function, the settling time, but this will have an insignificant effect on the solutions here).

\end{document}